\documentclass[a4paper, 12pt]{article}

\usepackage[norsk]{babel}
\usepackage[T1]{fontenc}
\usepackage[applemac]{inputenc}

\usepackage{amsmath,amssymb,amsfonts,amsthm}

\usepackage[pdftex]{hyperref}

\newcommand{\rmd}{\mathrm{d}}

\newcommand{\bbC}{\mathbb{C}}

\newcommand{\bbQ}{\mathbb{Q}}
\newcommand{\bbR}{\mathbb{R}}

\newcommand{\bbZ}{\mathbb{Z}}

\DeclareMathOperator{\SU}{SU}
\DeclareMathOperator{\SO}{SO}

\newcommand{\beq}{\begin{equation}}
\newcommand{\eeq}{\end{equation}}

\newcommand{\myemph}[1]
{ {\bf #1 }}

\newcommand{\mynewline}{\mbox{}\\}

\theoremstyle{definition}

\theoremstyle{remark}

\theoremstyle{plain}

\theoremstyle{definition}

\theoremstyle{remark}

\newcommand{\myquote}[3]
{
\vspace{2mm} 
\begin{flushright}
\begin{minipage}{9cm} 
 #3
\end{minipage}\\
\vspace{2mm} 
{\sc #1} \  {\it #2}
\end{flushright}
\vspace{2mm} 
}

\pdfpagewidth=\paperwidth
\pdfpageheight=\paperheight

\begin{document}

\title{
Matematikk ... Anvendelser
}

\author{Snorre H. Christiansen\footnote{
University of Oslo, Department of Mathematics, PO Box 1053 Blindern, NO-0316 Oslo, Norway.}\mbox{$\,$} \footnote{
I thank Knut Mørken for creating an arena for such discussions, in connection with the teaching reform InterAct at the  University of Oslo.
I am grateful to Christin Borge for innumerable conversations on these topics.}
}

\date{}

\maketitle

\abstract{Mathematics ... Applications. We consider a definition of mathematics as the art of thinking in terms of formalized systems, and the science of relations, structures and algorithms. We also touch upon the relation of mathematics to other sciences, in particular through modelling and scientific computing. We briefly discuss mathematics as a martial art and a key to paradigm changes. In Norwegian.}

\section{Innledningsvis}

\mynewline

\myquote{Charles Baudelaire}{Correspondances\footnote{A translation:\\
Nature is a temple where living pillars\\ 
Let sometimes emerge confused words; \\
Man crosses it through forests of symbols \\
Which watch him with intimate eyes.\\
--  Geoffrey Wagner, Selected Poems of Charles Baudelaire, NY: Grove Press, 1974.
}}{
La Nature est un temple où de vivants piliers\\
Laissent parfois sortir de confuses paroles ;\\
L'homme y passe à travers des forêts de symboles\\
Qui l'observent avec des regards familiers.
}

\newpage

\myquote{Hermann Weyl\footnote{This quote concludes \emph{''Why mathematics?'', you might ask},  by Michael Harris, page 966 in  \cite{Gow08}.}}{}{
With mathematics we stand precisely at the point of intersection of restraint and freedom that makes up the essence of man itself.
}

\myquote{Peter Lax}{Abel prize acceptance speech, 2005}{
Traditionally mathematics is divided into two kinds: pure and applied. The relation of the two is delicate. The great applied mathematician Joe Keller’s definition is: pure mathematics is a branch of applied mathematics. He meant that mathematics, beginning with Newton, was originally concerned with answering question[s] in physics, it is only later that the tools and concepts used were elaborated into theories that took on lives of their own. It was remarked by von Neumann that after a while abstract mathematics needs to be invigorated by the injections of new empirical material, like a new scientific theory, new experimental facts, or numerical studies.
}

\myquote{Michael Atiyah}{Preface to \cite{Arn00Frontiers}}{
For centuries, physics has had a close symbiotic relationship with mathematics, and Witten's forecast is that the 21st century will see this rise to new heights. But what of the future of biology? Many predict that understanding the brain will be the major challenge of the next century and, while it would be presumptuous of mathematicians to claim that they will solve the problem, it is not unreasonable to think that mathematics may have a useful part to play.
}

\newpage

\section{Hva skal barnet hete?}

Programnavnet \myemph{MIT} har svakheter. ''Informatikk'' er ikke sidestilt med ''Matematikk'' \emph{i dette programmet}.  Informasjonsoverfloden er en utfordring for oss, mens programmering er et middel. Matematikk er det overordnede prinsippet. ''Teknologi'' er for snevert i forhold til alle anvendelser vi ser for oss.  Det har også blitt poengtert at MIT navnet fungerer dårlig i rekrutteringsøyemed, fordi matematikken blir borte blant begrep som ikke henspiller på fag fra videregående. 

Vi behøver ikke sole oss i glansen av et teknologisk institutt i Massachusetts. Det er ikke bare i fredsmekling at ''Norway punches above its weight'', for å bruke Obamas uttrykk. Det er flere grunner til å kunne tro at våre kjære naboer i fremtiden vil si ''norrmännen har gjort det igen''. Norge har fostret flere matematikere av historisk kaliber (Abel, Lie, Skolem, Selberg ...).   I 2001 ble Ole-Johan Dahl og Kristen Nygaard tildelt Turing prisen for oppfinnelsen av objektorientert programmering. Matematisk Institutt har flere forskningsgrupper i verdenstoppen. Computers in Science Education bevegelsen arbeider for et helhetsperspektiv på utdannelse, som tar innover seg mulighetene datamaskiner gir, og har allerede fått anerkjennelse for sitt foregangsarbeid.

Som alternativer til MIT navnet vil jeg påstå at det ungdommen vil ha, og trenger, er MTV og MMA:

\begin{itemize}
\item \myemph{MTV} : i et samfunn dominert av overfladiskhet,  utgjør MatematikkTunge Vitenskaper et solid alternativ. Men MTV fungerer best som  fellesbetegnelse for flere programmer på Mat. Nat. fakultetet.

\item \myemph{MMA} : ''mixed martial arts'' er en underholdningsform bestående av voksne menn som sloss i oktogonale bur. For de av oss som er like fascinert av kampen som utspilles i andre mandalaer,  kan Matematikk Med Anvendelser være et vakrere og mer potent tidsfordriv.
\end{itemize}

\newpage

\section{Matematikk i tiden}

Matematikk er en kampkunst, på mer enn én måte. Matematikere redder liv og de tar liv, man trenger ikke gå lenger enn til \emph{The imitation game}, filmen om Alan Turing, for å illustrere dette. Men krigene er mangfoldige og ikke nødvendigvis åpenlyse.

Timothy Gowers rapporterte slik\footnote{På sin \href{https://gowers.wordpress.com}{Weblog}, August 25, 2014.} fra Emmanuel Candès' forelesning på ICM 2014:  ''He began by talking a little about the motivation for the kinds of problems he was going to discuss, which one could summarize as follows: his research is worth\-while \emph{because it helps save the lives of children}. ''  Det er her snakk om anvendelser av matematikk, nærmere bestemt ''compressed sensing'', i medisinsk billedbehandling. Universitetet i Oslo har plukket opp metoden, og man kan lese om dette i siste nummer av Apollon\footnote{No. 1/2015 side 28.}, under tittelen \emph{Medisinsk stråling kan reduseres til en sjettedel}.

Fra Arkimedes til von Neumann, har matematikere vært rådgivere for konger og presidenter, og bidratt til våpenutvikling og krigsstrategier. Per idag er NSA verdens største ansetter av matematikere. Vi kan være sikre på at også de har hørt om compressed sensing. Godfrey Harold Hardy, kanskje mest kjent for sin manglende evne til å forutsi tallteoriens potentialiteter\footnote{Se siste avsnitt i \cite{Har92}.}, oppgir at hans opprinnelige motivasjon for å drive med ''pure mathematics'', var at det ga han mulighet til å ''beat other boys''. Krystallkulen hans var muligens forkludret av et idyllisert syn på tallteoriens opprinnelse og historie. Fran\c cois Viète, oppfinneren av moderne algebra, var kong Henrik IV av Frankrikes kryptograf, lenge før RSA ble standarden for fortrolig kommunikasjon. I sitt bidrag til boken 
\cite{Arn00Frontiers}, utgitt i forbindelse med at år 2000 ble erklært matematikkens år, setter Vladimir Arnold det hele på spissen og hevder at all moderne matematikk er direkte koblet til krigføring. Roger Godement har inkludert en velinformert avhandling om matematikeres plass i det militærindustrielle komplekset i sin lærebokserie om Matematisk Analyse \cite{God98II}. 

Andre ''rene'' matematikere er heller ikke fremmede for krigens retorikk. Slik beskriver Claude Chevalley og Andr\'e Weil deres foregangsmann Hermann Weyl\footnote{Sitatet åpner forordet i \cite{Pes12}.}: ''A proteus who transforms himself ceaselessly in order to elude the grip of his adversary, not becoming himself again until after the final victory". Alle disse matematikerne så på matematikk som del av et intellektuelt og kulturelt landskap, og bidro til dette ikke bare med matematikk.  Vi noterer også at de trekker metaforene i retning av det filosofiske og guddommelige, i tråd med østerlandske prinsipper om kampkunst. 

Euklids verk Elementene var lenge kroneksemplet på matematisk tankegang og inspirerte tenkning langt ut over matematikkens grenser. Da Abraham Lincoln ble spurt om hvorfor han hadde ført opp på sin CV at han hadde lest Elementene, svarte han at han hadde ønsket å bli advokat, at som advokat var det nødvending å kunne argumentere, og at han ikke kunne tenke på noen bedre kilde for å lære argumentasjonsteknikk\footnote{Denne og andre anekdoter om matematikk og politikk er beskrevet i \cite{Kob14}.}. 
Spinoza presenterer sitt hovedverk ''Ethica'' som ''ordine geometrico demonstrata'' - altså bevist på geometrisk vis. Ikke alle matematikere er overbevist om at han lyktes i sistnevnte, men hans holistiske verdenssyn skinner kanskje klarest av alle den dag i dag. 

Disse eksemplene tyder på at innsikt i matematisk tankegang kan være til inspirasjon også når hovedanliggendet er å definere det gode liv og å styre staten. 

Matematikerens makt er skremmende. St. Augustin har uttalt\footnote{Dette sitatet er fra presentasjonen \emph{The Divine Madness --  Mathematics, Myths and Metaphors} av Morris W. Hirsch, 2008. Den inneholder flere tankevekkende utsagn.}    : ''Good Christians should avoid mathematicians and all impious soothsayers, taking care not to consort with those demons and deceptive spirits whose society will entrap them.''\footnote{Det latinske ordet ''mathematici'' \ i originalteksten hadde en videre betydning enn det vi idag forbinder med ''matematikere'' og inkluderte spesielt astronomer. Det er mulig det først og fremst var astrologer og numerologer St. Augustin ville advare mot.}
Hypatia, filosof, matematiker og dertil kvinne, ble i år 415 brutalt myrdet av en folkemengde oppildnet av munker, i et politisk oppgjør mellom rivaliserende kristne grupperinger med ulikt syn på fritenkningens verdi. At matematikk kan være så subversivt at livet står på spill, har både Kantorovitch (under Stalin) og Monge (under den Franske revolusjon) måttet ta  hensyn til\footnote{Se \cite{BarVil14} side 31 -- 32. Resten av denne dialogen, mellom en matematiker og en kunstner, i et buddhistisk tempel utenfor Paris, er også verdt å lese.}.

Finanskrisen anno 2008 skyldes ihvertfall delvis at finansielle instrumenter ble solgt til folk som ikke forsto matematikken som lå bak og dens begrensninger. Dermed hadde de ikke mulighet til å vurdere hvilken risiko de tok. Viljen var også svak. Regningen for etegildet ble jevnt fordelt, ikke bare på bordgjestene. 

Jean Leray, som på 30-tallet sto for milepæler innen teorien for ikkelineære partielle differensiallikninger og tolkningen av turbulens, valgte, da han ble holdt fange i Østerrike under annen verdenskrig, å presentere seg som spesialist i topologi, i frykt for at hans mer anvendbare kunnskaper skulle misbrukes av nazistene. Han levde opp til sitt ord, for under fangenskap organiserte han forelesninger og oppfant knipper og spektralsekvenser, begreper som senere har blitt sentrale for denne grenen av matematikken. 

Alexandre Grothendieck, spektralsekvensmagiker og en av det 20. århundrets matematiske visjonærer\footnote{\href{http://www.grothendieckcircle.org}{The Grothendieck Circle} er dedikert til å spre informasjon om hans liv og verk.}  gikk bort for ikke lenge siden. I kjølvannet av -68 revolusjonen, ble han så ubekvem med det vitenskapelige miljøets etikk at han etterhvert valgte å bryte helt med det. I denne nye perioden av hans liv var han med på å grunnlegge ''Survivre et vivre'', en aktivistgruppe med prinsipper ikke ulike Arne Næss' økosofi. 

\newpage

\section{Om matematikkens vesen}

Dette avsnittet er et forsøk på en tredelt definisjon av matematikk, etter følgende lest:
\begin{itemize}
\item Matematikk er kunsten å tenke ved hjelp av formaliserte systemer.
\item Som vitenskap har matematikk utviklet seg fra en refleksjon over størrelse og form, til læren om relasjoner, strukturer og algoritmer.
\item Matematiske språkformer, konsepter, metoder  og modeller er grunnleggende for de andre vitenskapene og utvikles i samspill med dem.
\end{itemize}

Forøvrig anbefales \cite{Gow02} som en lengre og mer konkret introduksjon til hva matematikk handler om. 


\subsection{Matematikk som kunst}

Matematikk har, i likhet med filosofi og naturvitenskap, sin opprinnelse i en undring over observerte lovmessigheter. De første tegnene til matematisk metode var en klargjøring av hva som utgjør gyldig argumentasjon og oppfinnelsen av diverse regnemetoder. 

Den første tråden fører frem til moderne logikk og mengdelære. Men heller enn at matematikk bør ses på som en gren av logikken, er matematisk logikk essentielt en anvendelse av matematisk tankegang på symbolstrenger, slik de fremstår i formelle språk. Matematisk logikk er en spesiell gren, fordi den inneholder miniatyrer av matematikken, og har avklart begreper som sannhet og bevis\footnote{For de fleste matematikere er dens resultater først og fremst negative, i den forstand at den har vist at noen formodentlig sanne påstander hverken kan bevises eller  motbevises.}. 

Det er bevisføring som gir matematikeren en trygg grunn å stå på. Men om bevis er gyldige dersom de kan oversettes til formelle bevis, så er det også et poeng at oversettelsen ville gjøre dem uforståelige for mennesker\footnote{''Human mathematics is a sort of dance around an unwritten formal text, which if written would be unreadable'' skriver David Ruelle i sitt bidrag til \cite{Arn00Frontiers}.}. Det vi jakter på i hverandres skrifter er intuisjoner og forståelse. Disse er vanskelige å definere, men vi kan assosiere dem med evnen til å reflektere over definisjoner (herunder utforske eksempler, moteksempler og varianter) og gleden av å se logiske sammenhenger, fra de \emph{såkalt} trivielle, til \emph{korrespondanser} på tvers av fagfelt\footnote{André Weil har skrevet til sin søster, filosofen Simone Weil, om analogiers rolle i hans matematikk \cite{Kri05}. På bakgrunn av slike intuisjoner formulerte Weil de formodingene Grothendieck senere brukte som rettesnor for sitt program for aritmetisk geometri.} \footnote{Denne publikasjonen ble gjenstand for noe så sjeldent som en åpen matematisk kontrovers, se \cite{Lan05Weil}.}.

Hvordan nå enn idéer gjemmer seg i teoremer og bevis\footnote{''We say that a proof is beautiful when it gives away the \emph{secret} of the theorem [...]'' (min uthevning) skriver Gian-Carlo Rota side 132 i \cite{Rot08}. }, er mendgespråket essentielt for måten matematikere kommuniserer på.

Den andre tråden fører frem til algoritmer og informatikk. Det fortettede symbolspråket som matematikere har innført for å utføre regnestykker,  har blitt videreutviklet til dagens programmeringsspråk. Koding og tungregning gjør det mulig å utføre regnestykker, både symbolske og numeriske, av en helt annen karakter enn tidligere.

Disse trådene er selvfølgelig flettet tett sammen. Bevisføring kan ses på som regning med idéer, og vi både kan og bør bevise egenskaper ved våre algoritmer. Formelle språk bør ses under ett, i tråd med Leibniz' drøm.

Matematikken utvider vår evne til rasjonell forståelse og forfiner vår intuisjon, på måter som gjør det mulig å utforske hittils ukjent territorium. Kavli-prisen\footnote{Se \href{http://www.kavliprize.org}{kavliprize.org} for mer informasjon.} nevner for eksempel uendelig små, store, eller komplekse fenomener, og da spesielt nano-, astro- og nevro- vitenskap, som våre store utfordringer.

Kvantefeltteori er et eksempel på formalisert system man kan regne med, uten at det hviler på et solid matematisk fundament. Helt frem til slutten av det 19. århundre, var situasjonen tilsvarende for kalkulus\footnote{Se \emph{The development of rigor in analysis} av Tom Archibald, side 117 i \cite{Gow08}.}. Å lage et tilfredsstillende fundament for kvantefeltteori er blant matematikeres store ambisjoner for det 21. århundre. Clay Mathematics Institute har også lovet ut dusør for oppdraget\footnote{Dette og de andre Millenium Problemene er beskrevet på \href{http://www.claymath.org}{claymath.org}.}.

\newpage

 \subsection{Matematikk som vitenskap}

En første innsikt er at matematiske objekter forstås best i måten de relaterer til hverandre, heller enn i sin ensomhet. Eksempler på relasjoner er funksjoner, komposisjonslover og ordener. De to siste kan forøvrig tolkes som eksempler på det første. Komposisjonslover tar to objekter av en gitt type og assosierer med dem, et nytt objekt av samme type. En ordensrelasjon uttrykker, for to gitte objekter av samme type, om den ene er større enn den andre på noen måte.

Komposisjonslover kan ses på som det grunnleggende begrepet i Algebra mens ordensrelasjoner har en tilsvarende plass i forhold til Analyse. Begge grenene, og gjerne i kombinasjon, er opptatt med å belyse Geometrien. Sistnevnte var både bemerkelsesverdig anvendbar og abstrakt.  Navigasjon til sjøs, ved hjelp av stjernene, er et eksempel på tidlig anvendelse. Spekulasjon over geometriske objekters natur var en av kildene til Platons idélære \cite{Nae01}.

En slik inndeling i Algebra, Analyse og Geometri fungerer godt på lavere grad, men er mindre relevant på høyere nivå, hvor grener hybridiseres, slik som i Algebraisk Topologi, eller knyttes nærmere anvendelsesområder, være seg bølger i nordsjøen eller kreftregisteret.

Matematikere liker å tenke på seg selv som mønsterfinnere\footnote{G. H. Hardys variant er ''A mathematician, like a painter or a poet, is a maker of patterns''.}, og man kan gjerne kalle slike regelmessigheter for strukturer.

Struktur kan også brukes mer restriktivt, slik tilfellet er i kategoriteori og objektorientert programmering. Mengder med en gitt struktur (f.eks. en komposisjonslov) utgjør gjerne en kategori og kategoriteori hjelper å klargjøre definisjoner og organisere flere typer matematiske strukturer (spesielt hjelper den med å overføre intuisjoner og tenkemåter mellom algebra og analyse). På tilsvarende måter er objektorientering nøkkelen til å organisere store kildekoder.

Funksjonsbegrepet er kanskje det mest elementære i matematikk, i den forstand at alt bygger på dens alkymi. En funksjon tilordner et nytt objekt  til hvert objekt av en viss type. Alle funksjoner er ikke konstruktivt definerbare, men for de som er det, er det som regel mange måter å regne dem ut på, noen vesentlig mer effektive enn andre. Etterhvert som mer komplekse fenomener beskrives matematisk, blir det nødvendig med nye avanserte algoritmer for å regne på dem. Regneferdigheter erstattes delvis med programmeringsferdigheter, men fremtiden krever også dypere matematisk forståelse av hvordan algoritmer fungerer og utvikles.

\newpage 

\subsection{Matematisk modellering}

Anvendelser og modeller er flertydige uttrykk. Man kan se på ren matematikk som en anvendelse av matematikk der estetiske prinsipper er rådende\footnote{Jf. bruken av ordet ''utelukkende'' i plakaten i inngangspartiet på N. H. Abels hus.}.  Man kan se på aritmetikk som en modell for det å telle og man kan se på algebra som en modell for aritmetikk. Logikere bruker ordet modell for noe som likner speilvendinger av disse påstandene. Vi velger å bruke ordet modellering i forbindelse med hvordan matematikken kommuniserer med andre vitenskaper. 

Man kan for eksempel modellere en golfball \cite{Arn12Golf} som et punkt som beveger seg i et tredimensjonalt rom, etter en gitt differensiallikning oppkalt etter Newton. Man kan velge om man vil ta hensyn til luftmotstand eller ei og man kan ha flere modeller for sistnevnte. Man kan  prøve å ta hensyn til usikkerhet knyttet til vindstyrke. Man kan kanskje erstatte punktet med en sphære og uttrykke at den beveger seg i en fluide beskrevet av partielle differensiallikninger oppkalt etter Euler. Man kan studere hvordan golfballen deformeres i det den blir slått ut, kanskje sprekker den, eller fortsetter å vibrere under sin flukt. I såfall kan man regne ut frekvensene, muligens høre dem. Man kan se på hvilken effekt de små gropene i overflaten til golfballen har på slagets presisjon, via turbulensen de skaper. Gitt modeller for alle usikkerhetsmomentene kan man studere sannsynlighetsfordelingen for nedslagspunktet. Gitt all denne kunnskapen kan man gjøre optimale valg på golfbanen.

Man snakker om hierarkier av modeller, med ulik kompleksitet og ulik grad av prediksjonskraft. Modeller kan være diskrete eller kontinuerlige, deterministiske eller stokastiske, og gjerne en blanding av de fire mulighetene dette gir\footnote{Dette er forøvrig en illustrasjon på forskjellen mellom $2 + 2$ og $2 \times 2$.}. Det kan være nødvendig å bruke flere modeller samtidig, for eksempel på forskjellige skalaer. Hvordan fenomener på forskjellige skalaer påvirker hverandre er idag et sentralt og uløst problem, belyst i Apollon fra flere vinkler.

Newtons nyvinning var imidlertid av en dypere karakter enn ballistikk anvendt på en peripatetisk sport. Han innførte kalkulus, et nytt matematisk språk, med dertil hørende teoremer og regnemetoder. Han postulerte nye naturlover som kunne uttrykkes i dette språket. Det gjorde det mulig å se eksperimenter på jorden i samme rammeverk som bevegelsene til himmellegemene, og gjøre forbausende presise prediksjoner i begge domener. Denne matematisk-fysiske teorien, godt hjulpet av mer reduksjonistisk modellering, ga ikke bare moment til den industrielle revolusjon, den forandret vårt verdensbilde. Kvantemekanikken, symbolisert ved Schrödinger likningen, utgjorde et tilsvarende paradigmeskifte i vår forståelse av verden, og banet vei for atomfysikken og dens teknologiske korollarer. Muligens krever  matematiseringen av biologien et nytt paradigmeskifte og muligens krever klimaforandringene et nytt Manhattan prosjekt.

\newpage 

\section{Aforismer om matematisk dannelse}

Titlene på de to neste avsnittene har jeg lånt fra en upublisert samling aforismer om India, skrevet av min far på 70-tallet. En tanke eller to, blant de bedre, kan jeg ha fra ham. 

\subsection{Innvielse}

\paragraph{Om fjellvett}
Vi liker å beklage oss over at studenter som begynner på universitetet hverken kan regne eller bevise. Men mange kan telle til fire og dette kan vi utnytte. Ta spesielt:
\begin{equation}
3 + 1 = 2 + 2.
\end{equation}
Leibniz fant utsagnet verdig en kommentar\footnote{\emph{Nouveaux essais sur l'entendement humain}, livre IV, chapitre VII.} og vi behøver ikke legge lista lavere. Definisjoner, teoremer og bevis handler om å si ifra om hvor man starter og hvor man er på vei, samt å ta passe lange skritt.

\paragraph{Om boksing}\footnote{I forbindelse med kontroversen over intuisjonismen skal Hilbert ha sagt: 
''Taking the principle of the excluded middle from the mathematician would be the same, say, as proscribing the telescope to the astronomer or to the boxer the use of his fists.''} Matematikeren har to never: den ekskluderte tredje og algoritmen. Noen foretrekker en venstre krok, andre en strak høyre, men ingen kommer langt med én hånd på ryggen.

\paragraph{Om å krysse Rubicon} 
Det er forøvrig en uting å telle helt til fire. Vi må lære våre studenter å telle:
\begin{equation}
\textrm{1, 2, mange,}
\end{equation}
ihvertfall hvis vi her erstatter mange med "$n$". Etter det, utfører man ikke ''for''-løkker for hånd, det har vi arméer av maskiner til.

\paragraph{Om å telle til 2}
Det er flere grunner til at man kan stoppe å telle når man har kommet til 2:
\begin{itemize}
\item irredusible polynomer over $\bbR$ har grad 1 eller 2.
\item minimalitet uttrykkes gjerne ved å studere deriverte av orden 1 og 2.
\end{itemize}
Kanskje har Newton også en finger med i spillet?

\paragraph{Om magisk notasjon}
Det følger av det over, at etter $n$ kommer $i$ og $\epsilon$. Den første er ikke imaginær, den er et meget konkret punkt i planet. Den andre er ikke infinitesimal, dens skjebne er å bli vilkårlig liten. Utfordringen vår er at disse likner til forveksling på $\sqrt{-1}$ og  $\rmd x$. Vi underviser ikke i Euler, men om.  Så kan man diskutere hva Robinson tilføyer.

\paragraph{Om tabuet rundt kropp} Vi har to kjøreregler:

-- Enhver teori krever to eksempler og ett moteksempel.

-- Fire eksempler krever en teori. 

\noindent La oss derfor ikke være blyge når vi snakker om $\bbQ, \bbR, \bbC, \bbR(X)$ med voksne mennesker, og ta $\bbZ_p$ i samme slengen.

\paragraph{Om de neste skritt}
Det skjer interessante ting også ved 3 og 4, for eksempel at $\SU(2)$ dekker over $\SO(3)$. Men for å ta steget fra det Euklidske plan til Lie grupper, bør man først kunne telle til 2, helt fra begynnelsen.

\subsection{Utvidelse}

\paragraph{Om kulturhistorie}\footnote{Les \cite{Stu00}.}
Du kan ikke telle til 5. Det har med norsk kulturhistorie å gjøre, og da tenker jeg ikke først og fremst på janteloven. Abels anerkjennelse kom sent.

\paragraph{Om skjønnhet} Matematikere har likevel tellet lenger enn til 4. Flate kristaller er det 17 av. Teorien er like vakker som flismønstrene i Alhambra og ikke urelatert. 

\paragraph{Om strandliv} Den ekstra friheten man får i den 5. dimensjonen gjør formenes mangfold mer håndterlig, ikke mindre. Det ble oppdaget av Smale i Rio på 60-tallet.

\paragraph{Om snarveier} Mye moderne matematikk handler om 24. Babylonerne ville vært henrykt. Men å telle så langt som til the Leech lattice krever ikke bare utholdenhet, det krever snarveier. Og arkeologisk finurlighet?

\paragraph{Om nirvana} Det var indiske brahminers innsikt at før $1$ kommer $0$. Hvor ellers, enn i deres fotspor, skal vi søke Dao?

\paragraph{Om uendelig} Fra Zenon og Archimedes til det algoritmiske og det tenkbare, har det potentielle og det aktualiserte antatt flere former. Men primtallene har vi ikke sett enden på, det har selv ikke Selberg.

\paragraph{Om selvet og mengden} \mbox{}

-- Det kraftigste argumentet er \emph{selvmotsigelsen}, for fra den følger alt.

-- Den dypeste innsikten stammer fra \emph{selvreferansen}, for i den aner vi vårt speilbilde.

Matematikkens gudommelige galskap er et system som tillater det siste uten å føre til det første. Stigen tar oss desto høyere. Tapre menn, deriblant Skolem, våget sin helse for å bringe oss denne flammen. 

\paragraph{Om dritt og emballasje}
Det har blitt sagt at gud skapte de naturlige tallene og at resten er menneskenes verk\footnote{Leopold Kronecker i en forelesning i 1886: ''Die ganzen Zahlen hat der liebe Gott gemacht, alles andere ist Menschenwerk''.}, bortsett fra fortegn, som er djevelens. Mange ble derfor misunnelige da von Neumann skapte det kumulative hierarki fra den tomme mengden. 

Dette systemet er \emph{no shit, only wrapping}. Om det gjør matematikken stueren, er det likevel ikke \emph{the real shit}. 

\paragraph{Om innsikt}
Det er godt å vite, men deilig å glemme. Matematikere setter dette i system, ved å identifisere strukturer, et sted i spennet mellom mengder og kategorier. 

\newpage

\section{Avslutningsvis}

\myquote{Apollon\footnote{Forskningsmagasinet Apollon, No. 1/2015, side 20.}}{Tema: beregninger, 2015}{
Penn og papir er ut. Avanserte beregninger på superraske datamaskiner er nå blitt helt nødvendig for å drive vitenskapen videre.
}

''Penn og papir er ut'' er et tvilsomt slagord når ''Beregnings-Universitetet'' lanseres. Steve Jobs oppga kalligrafi som det faget han lærte mest av på universitetet\footnote{Stanford commencement address, June 12, 2005.}. Man behøver ikke gå god for hvordan han har brukt sin lærdom, men vi må innse kraften i det som skjer når mennesker, deriblant matematikere, møter blanke ark med håndskrift. Det gjør det ikke mindre sant at beregninger spiller en økende rolle i alle vitenskaper. Kanskje vil spektralsekvenser komme godt med?

Vi står overfor gigantiske utfordringer knyttet til klima, energi og biologi. Matematikere kan spille nøkkelroller, som grønne teknologiutviklere, vitenskapelige innsiktsforedlere eller forbilder på kreativ og klar tankegang. 

Ikke alle er like entusiastiske til samfunnsutviklingen. Slik beskriver en visesanger fremtiden:

\myquote{Leonard Cohen}{The future}{
Things are going to slide, slide in all directions\\ 
Won't be nothing\\ 
Nothing you can measure anymore\\
The blizzard, the blizzard of the world\\ 
has crossed the threshold \\
and it has overturned \\
the order of the soul 
}

\vspace{5pt}

I kaoset får man lyst å synge med, dansende:

\myquote{Serge Gainsbourg}{L'Antillaise}{
Aux armes!\\
Et c\ae tera.
}

\newpage

\bibliography{../../Articles/Bibliography/newalexandria}{}
\bibliographystyle{plain}

\end{document}